\documentclass[conference]{IEEEtran}
\IEEEoverridecommandlockouts
\usepackage{cite}
\usepackage{hyperref}
\usepackage[english]{babel}

\usepackage{amsmath,amssymb,amsfonts}
\usepackage{algorithmic}
\usepackage{graphicx}
\usepackage{subfigure}
\usepackage{subfloat}
\usepackage{float}
\usepackage{textcomp}
\usepackage{xcolor}
\usepackage{textcomp}
\usepackage{array}
\usepackage{booktabs}
\usepackage{float}
\def\BibTeX{{\rm B\kern-.05em{\sc i\kern-.025em b}\kern-.08em
    T\kern-.1667em\lower.7ex\hbox{E}\kern-.125emX}}
    
\usepackage{etoolbox}
\makeatother
\makeatletter
\patchcmd{\@maketitle}
 {\addvspace{1\baselineskip}\egroup}
 {\addvspace{-1\baselineskip}\egroup}
 {}
 {}
\makeatother

\begin{document}

\title{Optimal offering strategy for an aggregator across multiple products of European day-ahead market
\thanks{This project is funded by the University of Adelaide industry-PhD grant scheme and Watts A/S, Denmark. This paper has been accepted for publishing by the conference IEEE PES ISGT Europe 2022.}
}

\author{\IEEEauthorblockN{Yogesh Pipada Sunil Kumar}
\IEEEauthorblockA{
\textit{University of Adelaide} \\
Adelaide, Australia \\
\href{mailto:yogeshpipada.sunilkumar@adelaide.edu.au}{yogeshpipada.sunilkumar},}
\and
\IEEEauthorblockN{S. Ali Pourmousavi}
\IEEEauthorblockA{
\textit{University of Adelaide}\\
Adelaide, Australia \\
\href{mailto:a.pourm@adelaide.edu.au}{a.pourm},}\and
\IEEEauthorblockN{Markus Wagner}
\IEEEauthorblockA{
\textit{University of Adelaide}\\
Adelaide, Australia \\
\href{mailto:markus.wagner@adelaide.edu.au}{markus.wagner@adelaide.edu.au}}
\and
\IEEEauthorblockN{Jon A. R. Liisberg}
\IEEEauthorblockA{
\textit{Watts A/S} \\
Svinninge, Denmark \\
\href{mailto:jon.lissberg@watts.dk}{jon.liisberg@watts.dk}}
}

\maketitle

\setlength\abovecaptionskip{-0.5\baselineskip}
\vspace*{-\baselineskip}

\begin{abstract}
Most literature surrounding optimal bidding strategies for aggregators in European day-ahead market (DAM) considers only hourly orders. While other order types (e.g., block orders) may better represent the temporal characteristics of certain sources of flexibility (e.g., behind-the-meter flexibility), the increased combinations from these orders make it hard to develop a tractable optimization formulation. Thus,  
our aim in this paper is to develop a tractable optimal offering strategy for flexibility aggregators in the European DAM (a.k.a. Elspot) considering these orders. Towards this, we employ a price-based mechanism of procuring flexibility and place hourly and regular block orders in the market. 
We develop two mixed-integer bi-linear programs: 1) a brute force formulation for validation and 2) a novel formulation based on logical constraints. 
To evaluate the performance of these formulations, we proposed a generic flexibility model for an aggregated cluster of prosumers that considers the prosumers' responsiveness, inter-temporal dependencies, and seasonal and diurnal variations.
The simulation results show that the proposed model significantly outperforms the brute force model in terms of computation speed. Also, we observed that using block orders has potential for profitability of an aggregator.
\end{abstract}


\section{Introduction}
\label{sec:Intro}
An aggregator acts as intermediary between prosumers (consumers who can produce) and the network operators (transmission and distribution). They manage a cluster of prosumers, aggregating behind-the-meter (BTM) flexibility from them to provide grid-scale ancillary services (voltage/frequency) to system operators. The operation of this entity provides the following benefits: 1) Increased utilization and ease of management of behind-the-meter (BTM) distributed energy resources (DER), 2) Providing cheap and fast flexibility services for system operators (transmission and distribution) improving reliability and aiding a larger integration of renewable energy sources in the generation mix, and 3) Boosting economic efficiency of the electricity markets~\cite{Burger2017}. Nonetheless, BTM flexibility aggregation has its own challenges such as uncertainty in BTM flexibility, aggregator's role in the market, market policy and optimization tools for simulating and evaluating different products, services and markets in which an aggregator may participate~\cite{7856228}. 

For making optimal ordering decisions in wholesale electricity markets, an aggregator must consider the uncertainty surrounding flexibility and wholesale prices as well as the nature of the flexible devices itself. Also, wholesale electricity markets are generally divided into day ahead markets (DAM) for scheduling generators for the next day and real-time/ancillary services markets to account for supply and demand mismatch. By participating in both markets, an aggregator can minimize costs compared with participation in DAM only~\cite{IRIA20191361}. This is because higher prices and penalties are applicable in real time markets, allowing the aggregator to hedge across these markets. Hence, most of the literature, comprehensively reviewed in~\cite{Okur2021}, apply stochastic optimization on multi-market aggregator problems to manage risk of aggregators under price and flexibility uncertainty. However, it would also be fruitful to consider the temporal characteristics of prosumer loads such as uninterruptible loads (washing machine, dryer, dishwasher) and batteries (EV and residential). The latter can degrade faster due to charge and discharge cycles at high power for short duration of time.

The European DAM, a.k.a. Elspot, allows participants to submit different order types such as hourly orders and block orders. An order is defined as a certain volume of energy bought/sold (offered/bid) at a certain price or price curve specified by the market participant for a certain duration. For hourly orders, volume is committed for a single hour only, whereas in \textbf{regular block order}, the same volume is committed for multiple consecutive hours. Block orders cannot be partially accepted; they must be accepted for the whole duration and hence, can represent the characteristics of certain sources of BTM flexibility (as stated above) and even unit commitment characteristics (UC) of thermal generators. Nonetheless, very little attention has been paid to them in the literature~\cite{9619919}. This is partially due to the combinatorial nature of these problems, making them computationally tedious to solve. 

The current state-of-the-art around block orders includes brute force models, where all the orders are explicitly provided to the optimization problem to solve \cite{Fleten2007}. To make it more tractable in some studies, the number of combinations/orders that can be placed is explicitly reduced \cite{Schledorn2021,Faria2011}. To the best of our knowledge, the only attempt at developing a non-brute force or non-explicit model for block orders is reported in~\cite{9619919}. They modeled different types of block orders by specifying a limit on the number of block orders that can be placed by a market participant for a thermal generator. Moreover, the above studies mainly focused on thermal generators, thus are unsuitable for flexibility aggregators. Therefore, there exists a gap in literature regarding block order based optimal bidding engines for flexibility aggregators.

To address this gap, we developed a novel mixed-integer bi-linear program formulation for a flexibility aggregator operating in Denmark to optimally offer BTM flexibility across all combinations of hourly and regular block orders in the Elspot market. We compared the performance of our model with the brute force model specified in~\cite{Fleten2007}. In addition, an appropriate flexibility model is developed for the aggregator considering the nature of the flexibility sources. Through extensive simulation studies, we show that our proposed formulation computationally outperforms the brute force model while achieving the same optimal results.

The rest of the paper is structured as follows. Section~\ref{sec:FlexProc} provides insights about the flexibility procurement for this aggregator. We then introduce our optimal problem formulations in Section~\ref{sec:Formulations}. The scenarios used for testing and validating the developed formulations are discussed in Section~\ref{sec:Simul}, and the corresponding results are presented in Section~\ref{sec:Results}. Finally, we draw conclusions and outline future work in Section~\ref{sec:Conclusion}.

\section{Flexibility Procurement}
\label{sec:FlexProc}
In this study, we assume that flexibility is procured from prosumers by using indirect load control strategy~\cite{Jordehi2019} wherein prosumers react to real-time price (RTP) signal sent by an aggregator. This is done to ensure prosumer autonomy is maintained. Flexibility procurement, thus, is divided into two components as described below. 

\subsection{Flexibility Model}
\label{subsec:Flex model}
Based on the products specified in~\ref{sec:Intro}, BTM sources can be categorized as block order BTM and non-block order BTM. The former includes sources such as EVs and residential batteries, where it is more beneficial to operate them for longer duration at same power (or energy) to lessen battery degradation. The latter contains thermal loads since they can only be flexible for a limited time before influencing the comfort of the prosumers. For the maximum benefit of the prosumers, PV is assumed non-curtailable, i.e, either used by the prosumers (may include battery storage) or sold to the grid.
The adopted flexibility model has two main components; a flexibility profile and a cross elasticity matrix.
Let \(\mathcal{T} = \{1, 2,\dots, T\}\) be the set of indices representing the operational hours, and \(T=\|\mathcal{T}\|\) be the number of operational hours. The flexibility profile \(F^{\text{max}} \in \mathbb{R}_-^{T \times 1}\) represents the maximum amount of flexibility that can be obtained at time interval $t$ in MWh, which depends on the weather conditions, diurnal factors and so on. Please note that we are not explicitly integrating these factors in our model since it is out of the scope of this paper. To account for inter-temporal relations (rebound effect) and prosumers' comfort, a cross-elasticity matrix \(A \in \mathbb{R}_-^{T \times T}\) is defined. This is a lower triangular matrix with diagonal elements equal to zero, which means that the flexibility used at time interval $t$ will affect the flexibility in the future intervals and not vice versa. \(F \in \mathbb{R}_-^{T \times 1}\) is a vector representing the actual flexibility obtained from the cluster of prosumers at time $t$. The superscripts $\text{hour}$ and $\text{block}$ represent the non-block order BTM and block order BTM, respectively. Since the aggregator can only make profit in the DAM by selling, all elements of the above variables are negative according to the Elspot's sign convention.
Equations~\eqref{eqref:Flexibility model A} - \eqref{eqref:Flexibility model C} represent the flexibility model used in the optimization. 
\begin{subequations}
\begin{flalign} 
\label{eqref:Flexibility model A}
    & F^{\text{max}, \text{hour}}  +  A^{\text{hour}}\, F^{\text{hour}}\leq F^{\text{hour}} \leq 0 &
\end{flalign}
\begin{flalign}
\label{eqref:Flexibility model B}
    & F^{\text{max}, \text{block}} + A^{\text{block}} \, F^{\text{block}} \leq F^{\text{block}} \leq 0 & 
\end{flalign}
\begin{flalign}
\label{eqref:Flexibility model C}
    & F = F^{\text{hour}} + F^{\text{block}} &
\end{flalign}    
\end{subequations}

\subsection{Prosumers' Responsiveness}
\label{subsec:Prosumer responsiveness} 
Relying on consumers' manual intervention to adjust their load according to price signal is not practical and efficient. Therefore, 
we assume that a home energy management system (HEMS) exists at the premises to manage communication, monitoring and control while fulfilling the prosumers preferences. Then, the prosumers' responsiveness refers to the price responsiveness of the rational prosumers, i.e., the price paid versus flexibility obtained. For a cluster of prosumers, this can be viewed as the supply curve of the aggregator that shows the cost to provide electricity to the wholesale market. 

\begin{figure}[ht]
    \centering
    \includegraphics[width=0.50\linewidth]{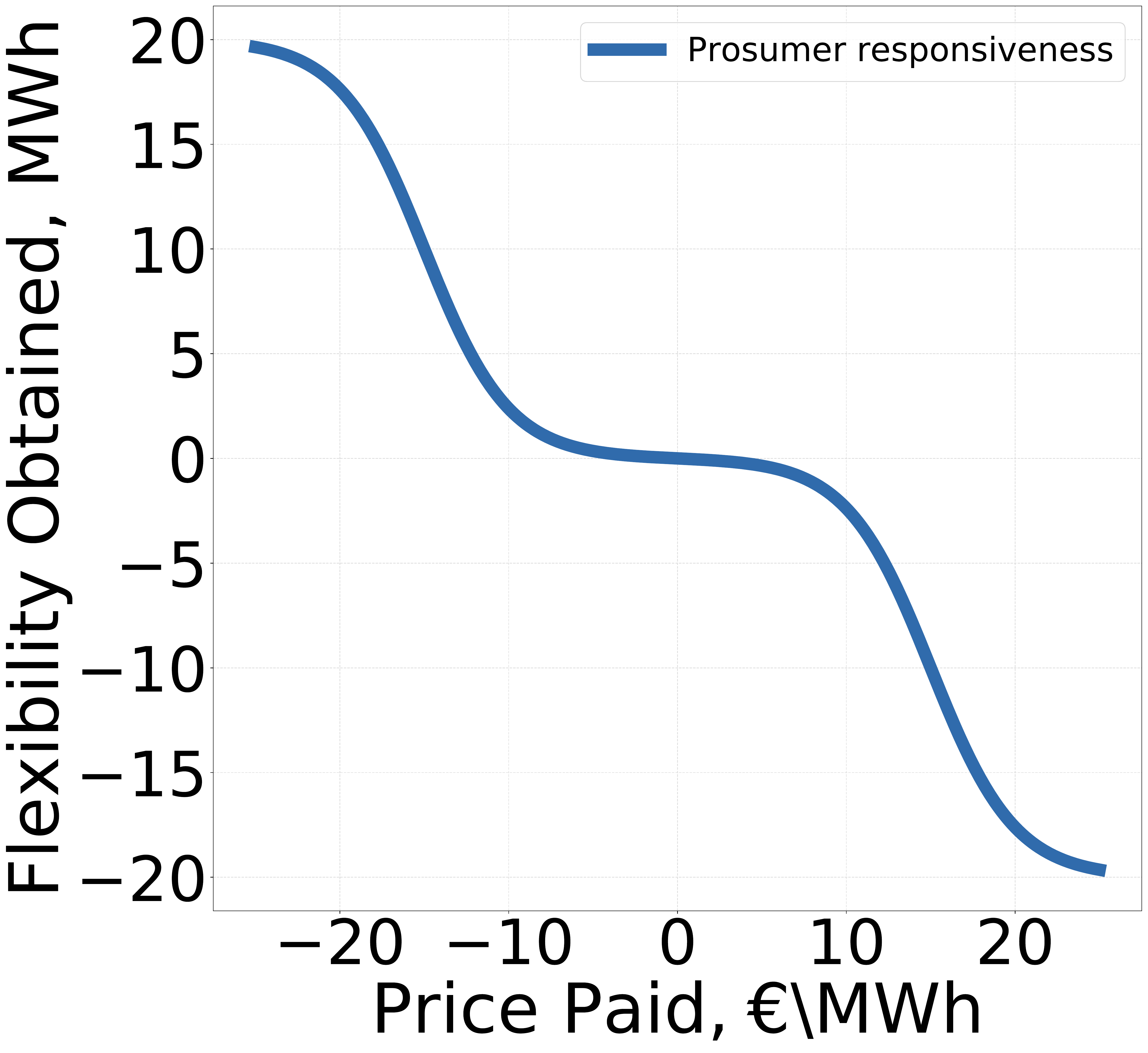}
    \caption{Expected Prosumer's responsiveness to different prices}
    \label{fig:Prosumer sigmoid}
\end{figure}

In general, prosumers responsiveness follows a saturation-like pattern, where little to no flexibility is available at low prices and no extra flexibility can be expected after reaching the maximum limit~\cite{8570785}. Thus, a sigmoid function can represent this pattern as shown in Fig.~\ref{fig:Prosumer sigmoid}. 

Here, we use the same price responsiveness curve for a prosumer in all operational hours to simplify the problem, but the curve can be modified for different hours if needed. Let \(\lambda^{\text{paid}} \in \mathbb{R}^ {T \times 1} \) represent the price paid by the aggregator at each operational hour, \(a \text{ and } b\) are constants associated with prosumer responsiveness and \(f^{\text{max}} \text{ and } \lambda^{\text{max}}\) denote the knee point flexibility and the corresponding price paid for it, respectively, from a cluster of prosumers for all operational hours. 
\begin{flalign}
    \label{eqref:Price paid vs Flexibilility obtained sigmoid}
    & F = -\frac{f^{\text{max}}}{1 + e^{\circ\left( a + b \, \lambda^{\text{paid}}\right)}} &
\end{flalign}
Equation~\eqref{eqref:Price paid vs Flexibilility obtained sigmoid}\footnote{Please note that \(\circ\) symbol represents element-wise or Hadamard operation.} is a non-convex function when double/single sided. To make the overall optimization problem tractable, this function is linearized by piece-wise linearization (PWL). Let \(\mathcal{P} = \{1, 2,\dots, P\}\) be the set of indices representing each piece in the PWL model, where \(P=\|\mathcal{P}\|\) is the number of pieces, \(\lambda^{\text{PWL}} \in \mathbb{R}^{P+1 \times 1}\) is the list of price breakpoints, \(F^{\text{PWL}} \in \mathbb{R}^{P+1 \times 1}\) and \(s^{\text{PWL}} \in \mathbb{R}^{P+1 \times 1}\) are the corresponding lists of flexibility and slopes obtained from Eq.~\eqref{eqref:Price paid vs Flexibilility obtained sigmoid}. To apply a \(P-\text{piece}\) linearization, the following linear constraints are introduced:
\begin{subequations}
\begin{flalign}
    \label{eqref: PWL A}
    & \sum_{p \in \mathcal{P}} b_{t, p}^{\text{PWL}} = 1 
    & 
\end{flalign}
\begin{flalign}
    \label{eqref: PWL B}
    & u_{t,p}^{\text{PWL}} \leq b_{t, p}^{\text{PWL}} \cdot \left(\lambda^{\text{PWL}}_{p+1} - \lambda^{\text{PWL}}_{p}\right) \quad \forall \: t \in \mathcal{T}, p \in \mathcal{P} & 
\end{flalign}
\begin{flalign}
    \label{eqref: PWL C}
    &\lambda_{t}^{\text{paid}} = \sum_{p \in \mathcal{P}} \left( u_{t,p}^{\text{PWL}} + b_{t,p}^{\text{PWL}} \, \lambda^{\text{PWL}}_{p}\right) \quad \forall \: t \in \mathcal{T} & 
\end{flalign}
\begin{flalign}
    \label{eqref: PWL D}
    & F_t = \sum_{p \in \mathcal{P}} \left[b_{t,p}^{\text{PWL}} \cdot \Big(F^{\text{PWL}}_{p} - F^{\text{PWL}}_{1}\right) + s_p ^{\text{PWL}} \, u_{t,p} ^{\text{PWL}} \Big]  &  \nonumber\\ 
    & \hspace{6.5cm} \forall \: t \in \mathcal{T} & 
\end{flalign}
\begin{flalign}
    \label{eqref:PWL E}
    &-f^{\text{max}} \leq F_t & 
\end{flalign}
\begin{flalign}
    \label{eqref:PWL F}
    &0 \leq \lambda^{\text{paid}}_t \leq \lambda^{\text{max}}&
\end{flalign}
\end{subequations}
\noindent where \(b^{\text{PWL}} \in \{0, 1\}^{T \times P}\) is a matrix of binaries, which corresponds to the ``piece'' that \(\lambda_{t}^{\text{paid}}\) is in at time \(t\), \(u^{\text{PWL}} \in \mathbb{R}^{T \times P}\) is a matrix of auxiliary variables used to calculate the value of \(\lambda_{t}^{\text{paid}}\) at time \(t\).  

\section{Optimization problem formulations}
\label{sec:Formulations}
This section describes the two optimal offering strategies developed for maximizing aggregator's profit, while satisfying the market rules and flexibility constraints. An offer is defined by the volume offered (price-independent case), a start time of delivery and the duration of delivery. The optimization model must consider all possible combinations of single-hour order and regular block orders to choose the most profitable ones. These combinations can be modeled by the use of binary variables, where different types of order cannot be placed at a given operational hour. The profit of the aggregator contains two components: revenue from the market and cost of flexibility procurement, both of which are bi-linear, i.e., (\(\text{price} \times \text{volume}\)). This requires developing a mixed-integer bi-linear program (MIBLP), which are NP-hard problems, and that can be solved using modern solvers like Gurobi~\cite{Gurobi21}.

In this study, to improve the tractability of the problem, the following assumptions are made: 1) the aggregator is a price taker, and 2) the formulation is deterministic, i.e., the prosumer flexibility and market clearing prices are exactly known at the start of the optimization. We recognize that the flexibility model is highly uncertain and that modeling this uncertainty is critical for developing aggregator's business model. However, the main purpose of this study being the introduction of a new way of modeling block orders, we postpone uncertainty modeling to the future. 

Based on the above assumptions, the developed optimization models are explained below. 

\subsection{Brute Force Model}
\label{subsec:Brute Force}
In this model, all possible combinations of single and block orders are pre-defined for the optimization problem along with the flexibility constraints and the spot price, \(\lambda^{\text{spot}} \in \mathbb{R}^{T \times 1}\), from the Elspot market. The optimization problem maximizes aggregtor's profit by selecting the most profitable combination of orders. This approach has been implemented in~\cite{Faria2011, Schledorn2021, 9619919}.
\renewcommand{\footnotesize}{\scriptsize}
Let \(\mathcal{H} \in \{0, 1\}^ {T \times H}, \mathcal{B} \in \{0, 1\}^{T \times B}\) be the set of binary vectors representing combinations of hourly and block orders respectively\footnote{$ \mathcal{H} = \Big\{ \big\{1, 0 , 0, \dots,0, 0, 0\big\},\big\{0, 1 , 0\dots,$\\
\hspace*{1.45cm}$0, 0, 0\big\}, \dots ,\big\{0, 0 , 0, \dots,0, 1, 0\big\},\big\{0, 0, 0, \dots, 0, 0, 1\big\}\Big\}$\\
\hspace*{0.4cm}$\mathcal{B} = \Big\{ \big\{1, 1 , 1, 0 , \dots,0, 0, 0\big\},\big\{0, 1 , 1, 1, 0, \dots,$\\
\hspace*{1.41cm}$0, 0, 0\big\},   \dots ,\big\{0, 1 , 1, 1, 1  \dots,1, 1, 1\big\},\big\{1, 1, 1, 1, 1 \dots, 1, 1, 1\big\}\Big\}$
}. \(H\) and \(B\) are the number of combinations of hourly and block orders, respectively. \(V^{\text{hour}} \in \mathbb{R}_-^{H \times 1}\) is the volume offered for each hourly order combination in \(\mathcal{H}\), while \(V^{\text{block}} \in \mathbb{R}_-^{B \times 1}\) is the volume offered for each block order combination in \(\mathcal{B}\).
\(b^{\text{hour}} \in \{0, 1\}^{H \times 1}\) and \(b^{\text{block}} \in \{0, 1\}^{B \times 1}\) are two sets of binary vectors encoding which combination of hourly and block orders are active, respectively. Thus, the formulation is as follows,

\begin{subequations}
\begin{flalign}
    \label{eqref:Objective 1}
    & \max_{\Psi_1} \sum_{t \in \mathcal{T}} -F_t \cdot \left(\lambda^{\text{spot}}_t - \lambda^{\text{paid}}_t\right) &
\end{flalign}
\begin{flalign}
    & \textbf{subject to:} \text{ Eqs.~}\eqref{eqref:Flexibility model A}-\eqref{eqref:Flexibility model C}, \eqref{eqref: PWL A} - \eqref{eqref:PWL F} & \nonumber  
\end{flalign}
\begin{flalign}
    \label{eqref:Hourly order Flexibility}
    & F_t ^ {\text{hour}} = \sum_{h=1}^{\|\mathcal{H}\|} V_h^{\text{hour}}\,\mathcal{H}_{t,h}  \quad \forall \: t \in \mathcal{T}& 
\end{flalign}
\begin{flalign}
    \label{eqref:Block Order Flexibility}
    & F_t ^ {\text{block}} = \sum_{b=1}^{\|\mathcal{B}\|} V_h^{\text{block}}\, \mathcal{B}_{t,b}  \quad \forall \: t \in \mathcal{T}& 
\end{flalign}
\begin{flalign}
    \label{eqref:One order active 1}
    &\sum_{h=1}^{\|\mathcal{H}\|} b_h ^{\text{hour}} \, \mathcal{H}_{t,h} + \sum_{b=1}^{\|\mathcal{B}\|} b_b ^{\text{block}} \, \mathcal{B}_{t,b} \leq 1 \quad \forall \: t \in \mathcal{T}&
\end{flalign}
\end{subequations}

\noindent where \(\Psi_1 = \{F, F^{\text{hour}}, F^{\text{block}}, V^{\text{hour}}, V^{\text{block}}, b^{\text{hour}}, b^{\text{block}} \linebreak[0] ,\lambda^{\text{paid}}\}\) is the set of decision variables to be optimized, Eq.~\eqref{eqref:Objective 1} is the profit of the aggregator, Eqs.~\eqref{eqref:Hourly order Flexibility} and \eqref{eqref:Block Order Flexibility} obtain the vector of flexibility ordered at each operational interval for hourly and block orders, respectively. Finally, Eq.~\eqref{eqref:One order active 1} ensures that at a given operational interval, the maximum of one order is active out of all the combinations in \(\mathcal{H} \text{ and } \mathcal{B}\). 

\subsection{Proposed Formulation}
\label{subsec:AM2}

\begin{figure}[ht]
    \centering
    \includegraphics[clip, trim = {5.1cm, 4.1cm, 4.1cm, 4.1cm}, width=0.75
    \linewidth]{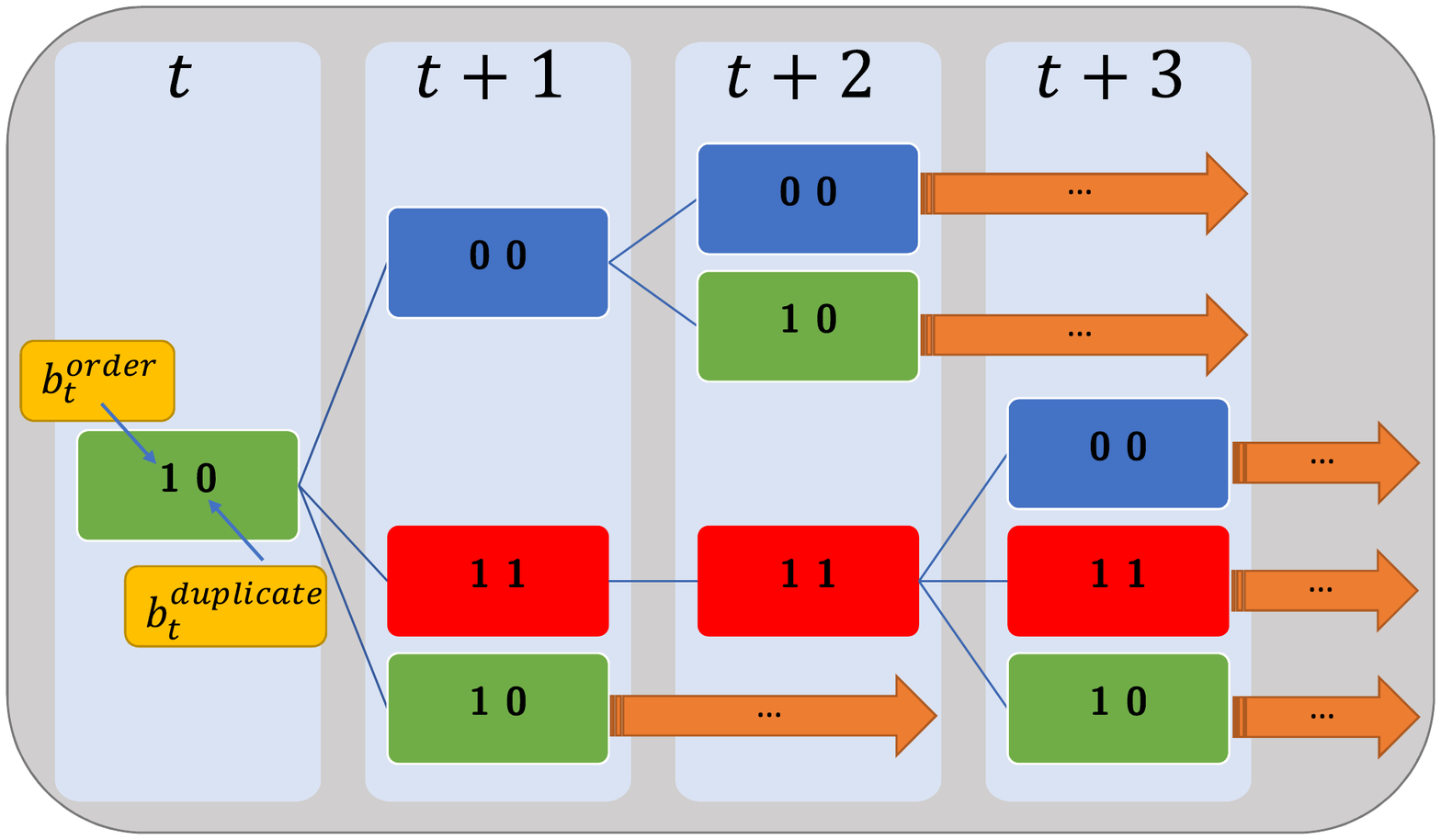}
    \caption{The implemented decision tree of the proposed formulation}
    \label{fig:Decision Tree}
\end{figure}

The new formulation is based on the fact that the only difference between regular block orders and hourly orders is the equal volume for a minimum duration of three intervals. Thus, it is sufficient to track duplication of orders consecutively to define a block order. 

\(b^{\text{order}}_t\) signifies whether an hourly or block order is placed at time \(t\) or not. Let \(b^{\text{duplicate}} \in \{0, 1\}^{(T+1) \times 1}\) be a binary vector indicating whether the previous order is duplicated or not, and \(u^{\text{active}} \in [0, 1]^{(T+1) \times 1}\) be an auxiliary vector used to enforce the minimum duration of three hours for a block order based on the Elspot rules. The above variables are indexed using the index set \(\{0\} \cup \mathcal{T}\). The values of \(b^{\text{order}}_0, b^{\text{duplicate}}_0\), \(u^{\text{active}}_0\) and \(u^{\text{active}}_1\) are redundant and set to zero, and are used to generalize the constraints.  
\(u^{\text{hour}}, u^{\text{block}} \in [0, 1]^{T \times 1}\) are two vectors of binary variables encoding hourly and block orders status at all operational intervals, respectively. The start of a block order is tracked using a binary variable, \(u^{\text{start}} \in [0, 1]^{T \times 1}\), which is set to one only if there is a transition from a non-block order to a block order. It is noted that the above three binary variables are relaxed to real variables. This is because the logical constraints mentioned in the formulation below are reformulated into MILP constraints as shown in~\cite{WolffTM14}, which ensures that the solution of these variables are always at the boundary. 

\(b^{\text{vol},+} \in \{0, 1\}^{T \times 1}\) and \(b^{\text{vol},-} \in \{0, 1\}^{T \times 1}\) are binary variables used to enforce \(F_{t-1} = F_{t}\), if both are equal to 1, by using \textit{Big M} formulation. Figure~\ref{fig:Decision Tree} shows the decision tree implemented via the logical constraints in the proposed formulation. The actual interpretations of each state are explained in Table~\ref{tab:Block order}.
\begin{table}[ht]
\small
    \centering
        \caption{Proposed model binary variable interpretation}
        \vspace{0.1cm}
    \label{tab:Block order}
    \begin{tabular}{lll}
    \toprule
    \(b^{\text{order}}_t\) & \(b^{\text{duplicate}}_t\) & \textbf{Interpretation} \\
        \midrule
         0 & 0 & No order is placed\\
         1 & 0 & Some order is placed\\
         1 & 1 & Previous order is duplicated\\
         0 & 1 & Illegal state\\
        \bottomrule
    \end{tabular}
\end{table}

Thus, the optimization problem can be expressed as follows, 
\begin{subequations} 
\begin{flalign}
&\max_{\Psi_2} \sum_{t \in \mathcal{T}} -F_t \cdot \left(\lambda^{\text{spot}}_t - \lambda^{\text{paid}}_t\right)& 
\end{flalign}
\begin{flalign}
&\textbf{subject to:} \text{ Eqs.~}\eqref{eqref:Flexibility model A}-\eqref{eqref:Flexibility model C}, \eqref{eqref: PWL A} - \eqref{eqref:PWL F}&  \nonumber  
\end{flalign}
\begin{flalign}
    \label{eqref:initialize 3}
    &b_0^{\text{order}}, b_0^{\text{duplicate}}, u_0^{\text{active}}, u_1^{\text{active}} = 0 &
\end{flalign}
\begin{flalign}
    \label{eqref:Block order duration encoding 3}
    &u_{t}^{\text{active}} = b_{t - 1}^{\text{duplicate}} \wedge \neg b_{t - 2}^{\text{duplicate}} \, \forall \: t \in \mathcal{T} - \{1\} & 
\end{flalign}
\begin{flalign}
    \label{eqref:Illegal sequence encoding 3A}
  &  \neg b_{t - 1}^{\text{order}} \wedge b_{t}^{\text{duplicate}} = 0 \quad \forall \: t \in \mathcal{T} &
\end{flalign}
\begin{flalign}
  \label{eqref:Illegal sequence encoding 3B}
  &  \neg b_{t}^{\text{order}} \wedge b_{t}^{\text{duplicate}} = 0 \quad \forall \: t \in \mathcal{T} &
\end{flalign}
\begin{flalign}
    \label{eqref:Enforcing block order constraints 3} 
    &u_{t}^{\text{active}} \leq b_{t}^{\text{duplicate}} \leq u_{t}^{\text{active}} + 1 \quad \forall \: t \in \mathcal{T}&
\end{flalign}
\begin{flalign}
    \label{eqref:Block order start encoding 3A}
    &u_{t}^{\text{start}} = \neg b^{\text{duplicate}}_t \wedge b^{\text{duplicate}}_{t+1} \, \forall \: t \in \mathcal{T} & 
\end{flalign}
\begin{flalign}
    \label{eqref:Block order start encoding 3B}
    &u^{\text{start}}_t = 0 \quad \forall \: t \in \{T - 1, T\} &
\end{flalign}
\begin{flalign}
    \label{eqref:Block order status encoding 3}
    & u^{\text{block}}_t = \left(b^{\text{order}}_t \wedge b^{\text{duplicate}}_{t} \right) \vee u^{\text{start}}_t \quad \forall \: t \in \mathcal{T} \, \backslash \{T\}& 
\end{flalign}
\begin{flalign}
    \label{eqref:Single order status encoding 3}
    & u_{t}^{\text{hour}} + u_{t}^{\text{block}} = b^{\text{order}}_t \quad \forall \: t \in \mathcal{T} &
\end{flalign}
\begin{flalign}
    \label{eqref:Equal Volume encoding 3}
   & b^{\text{duplicate}}_{t} = b^{\text{vol},+}_{t} \wedge b^{\text{vol},-}_{t} \quad \forall \: t \in \mathcal{T}&
\end{flalign}
\end{subequations}
\noindent where, $\Psi_2 = \{F, F^{\text{hour}}, F^{\text{block}}, b^{\text{order}}, b^{\text{duplicate}}, b^{\text{vol},+} , b^{\text{vol},-},  \linebreak[0] \lambda^{\text{paid}}, u^{\text{active}}, u^{\text{block}}, u^{\text{hour}}, u^{\text{start}}\}$ is the set of decision variables in the optimization problem. Equation~\eqref{eqref:initialize 3} initializes all the binary and relaxed binary variables to zero at \(t = 0\). Equation~\eqref{eqref:Block order duration encoding 3} is used to check if a block order is within the minimum duration constraint or not. Equation~\eqref{eqref:Enforcing block order constraints 3} enforces an active block order if inside minimum duration. Equation~\eqref{eqref:Equal Volume encoding 3} is used to evaluate and enforce the equal volume constraint for the block orders. Equations~\eqref{eqref:Illegal sequence encoding 3A}-\eqref{eqref:Illegal sequence encoding 3B} block the illegal sequences mentioned above. Equation~\eqref{eqref:Block order start encoding 3A} encodes the \(u^{\text{start}}_t\), while Eq.~\eqref{eqref:Block order start encoding 3B} prevents block orders starting at times when minimum duration requirement cannot be satisfied. Also note that this formulation can model \textit{profile block orders}~\cite{DAMOrder}) as well, by dropping Eq.~\eqref{eqref:Equal Volume encoding 3} and its associated \textit{Big M} formulation.

\begin{figure*}[h]
        \centering
            \subfigure[Hourly order flexibility scenarios]
            {
                \label{subfig:Single order flex}
                \includegraphics[width=0.23\textwidth]{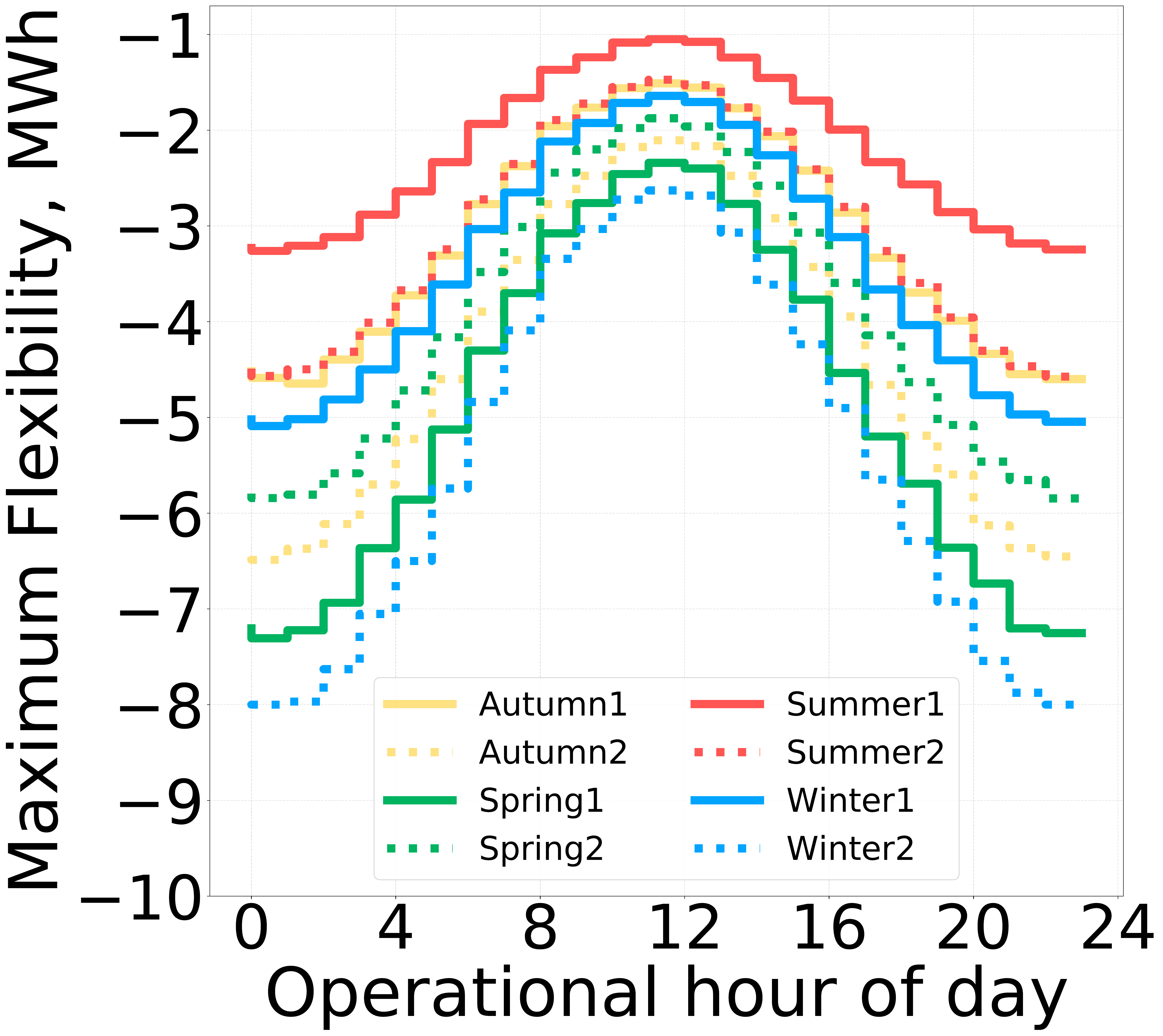} 
            } 
            \subfigure[Block order flexibility scenarios] 
            {
                \label{subfig:Block order flex}
                \includegraphics[width=0.23\textwidth]{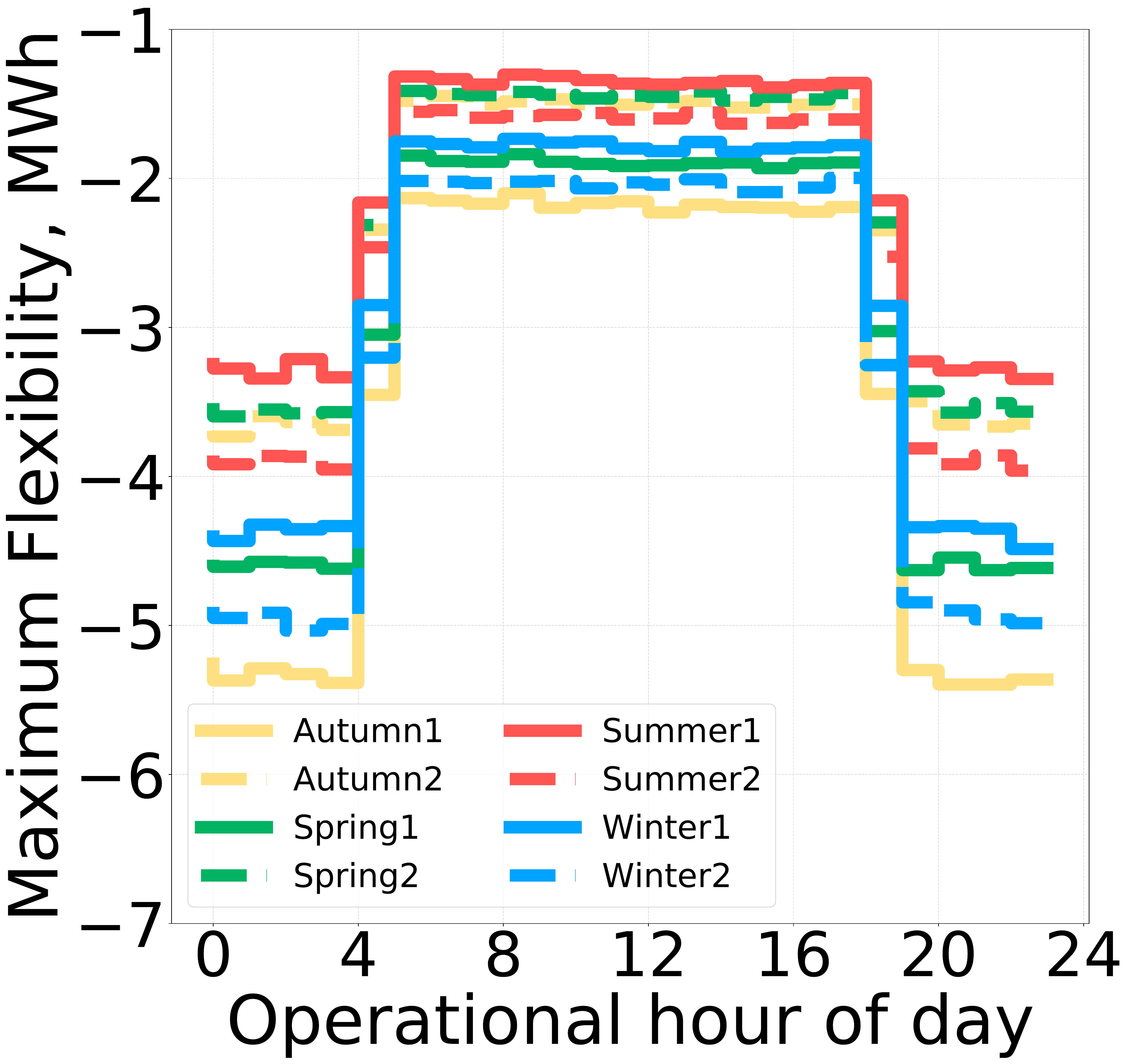} 
            } 
            \subfigure[Elspot market price scenarios] 
            {
                \label{subfig:Elspot scenarios}
                \includegraphics[width=0.23\textwidth]{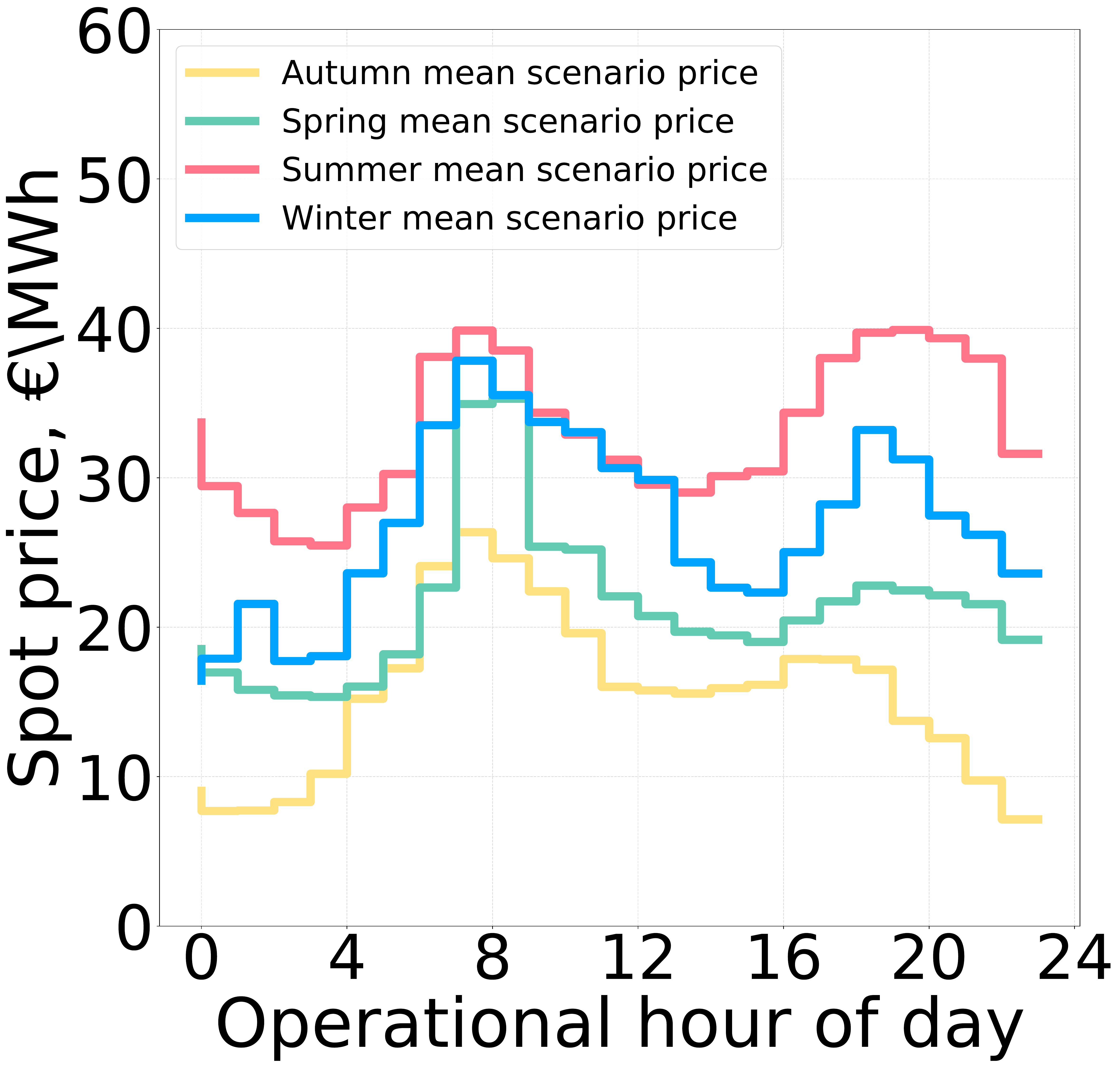}
            }
            \subfigure[Prosumers' responsiveness scenarios with \(a=6\), \(b=-0.4\)] 
            {
                \label{subfig:Prosumer responsiveness PWL}
                \includegraphics[width=0.23\textwidth]{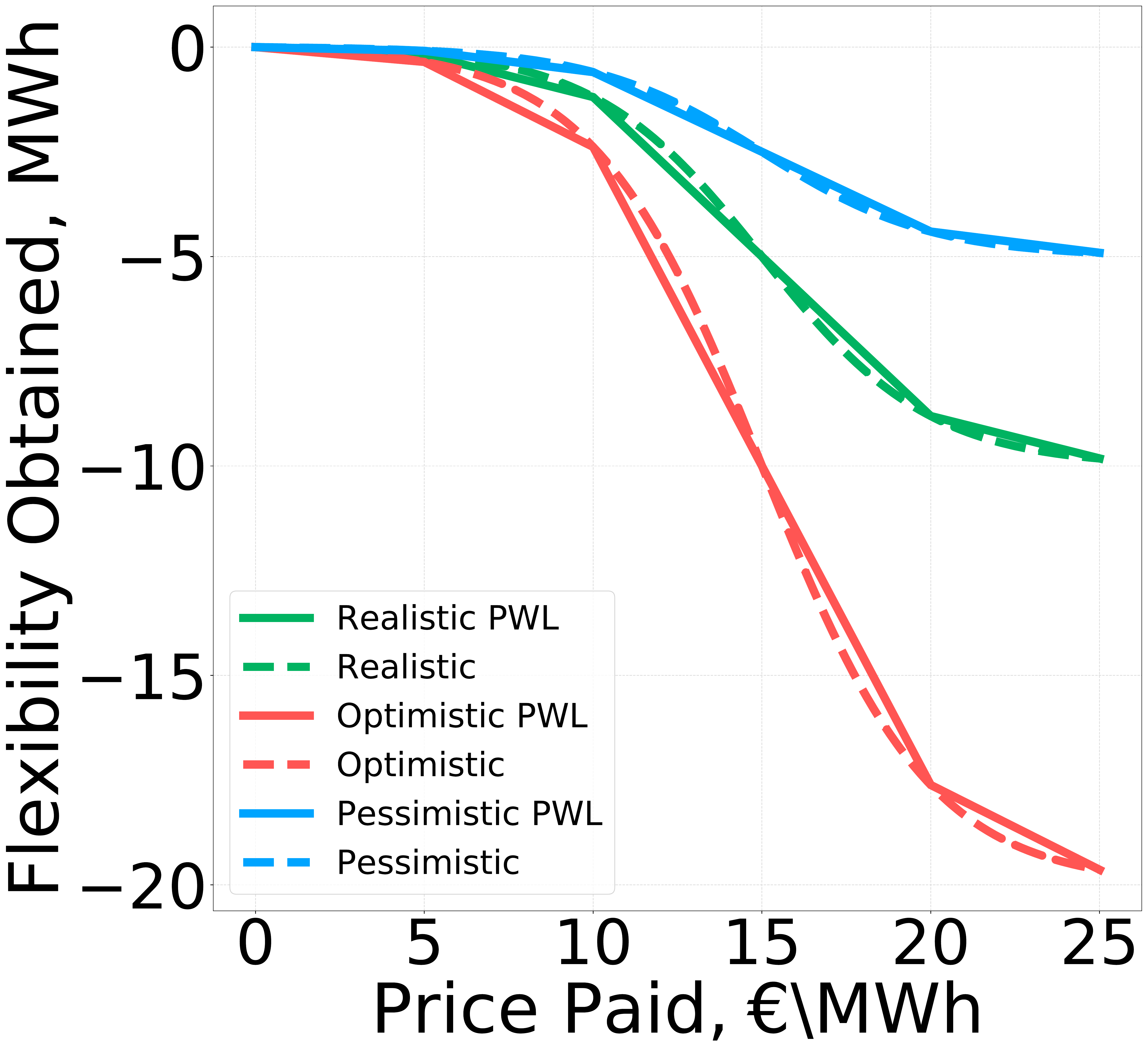}
            }
        \caption{Simulation scenarios}
        \label{fig:Scenarios}
\end{figure*}

\section{Simulation Setup}
\label{sec:Simul}
To compare our three models, we use multiple scenarios of flexibility and prosumer responsiveness based on real Elspot price data. 
Denmark was chosen as the region of operation of the aggregator. 

\subsection{Flexibility scenarios}
\label{subsec:Simul flexibility}
As discussed in Section~\ref{subsec:Flex model}, one flexibility profile and cross-elasticity model is needed for block orders and non-block orders per scenario. Prosumers' flexibility mainly depends on factors such as seasonal and geographical variations, diurnal factors (time of day), and occupant behavior, of which seasonal variations have the highest effect as reported in~\cite{Balint2019}. This study was used as a basis to design cases for non-block order profiles for different seasons and diurnal factors as shown in \ref{subfig:Single order flex}. In case of block order flexibility, seasonal variations were not considered and the names in Fig.~\ref{subfig:Block order flex} are the hourly order profiles with which these block order profiles are paired. Since these include EV, the flexibility would be maximum when prosumers are at home (7 PM to 4 AM). Hence, eight profiles and cross-elasticity matrices (CEM) were generated as scenarios, two for each season (e.g., \(Autumn1, Autumn2\)). The profiles and CEM were generated using log-normal distribution to ensure one-sided profiles and matrix elements (selling only) as explained in Section~\ref{subsec:Flex model}. Please note that proper scenarios can be generated using historical values based on the prosumers' reaction at different times of a day and seasons but that is outside the scope of this study.

Additionally, for each profile and cross-elasticity matrix pair, two Elspot price scenarios (\textit{Energy Price 1, Energy Price 2}) were considered while simulating the cases. These prices were chosen from the DK1 region of Denmark~\cite{DKKprice}, such that the prices belonged to the months of the corresponding season of the profile in the year 2021. Figure~\ref{subfig:Elspot scenarios} shows the average Elspot prices for each season (four scenarios per season). Thus, 16 scenarios combining prosumer flexibility and Elspot prices were generated in total. \(T=24, B= 253, H=24, M=40, P=5\) are the parameters used for simulation. 
\vspace{-0.3em}
\subsection{Prosumers' responsiveness scenarios}
\label{subsec:Simul resposiveness}
As discussed in Section~\ref{subsec:Prosumer responsiveness}, Eq.~\eqref{eqref:Price paid vs Flexibilility obtained sigmoid} is used to model prosumer responsiveness where \(a\) and \(b\) can be viewed as shaping parameters while \(f\) limits the maximum flexibility. The average electricity spot prices between 2010 to 2020 for DK1 region for Elspot market was found to be \texteuro 35.24/MWh. The function specified in Eq.~\eqref{eqref:Price paid vs Flexibilility obtained sigmoid} represents the cost curve of the aggregator to activate flexibility and the ``knee point'' of this is at the price of \texteuro 25.00/MWh. This is in the aggregators' interest as it would be more likely to generate profit from the wholesale market. Also, it provides a decent incentive for prosumers to participate in the flexibility aggregation program. Additionally, based on the maximum available flexibility, three scenarios of prosumers' clusters are simulated; \textit{Optimistic, Realistic and Pessimistic} with \(f = 20, 10, 5\), respectively. Finally, they are PWL approximated using the \(\lambda^{\text{PWL}} = [0, 5, 10, 15, 20, 25]\), as shown in Fig.~\ref{subfig:Prosumer responsiveness PWL}.  

Thus, we consider a total of 48 scenarios by combining the 16 scenarios explained in Section~\ref{subsec:Simul flexibility} and the three prosumers' responsive scenarios discussed here. 

\begin{figure*}[h]
        \centering
            \subfigure[Objective value comparison for acceptable scenarios]
            {
                \label{subfig:Obj comparison}
                \includegraphics[width=0.23\textwidth]{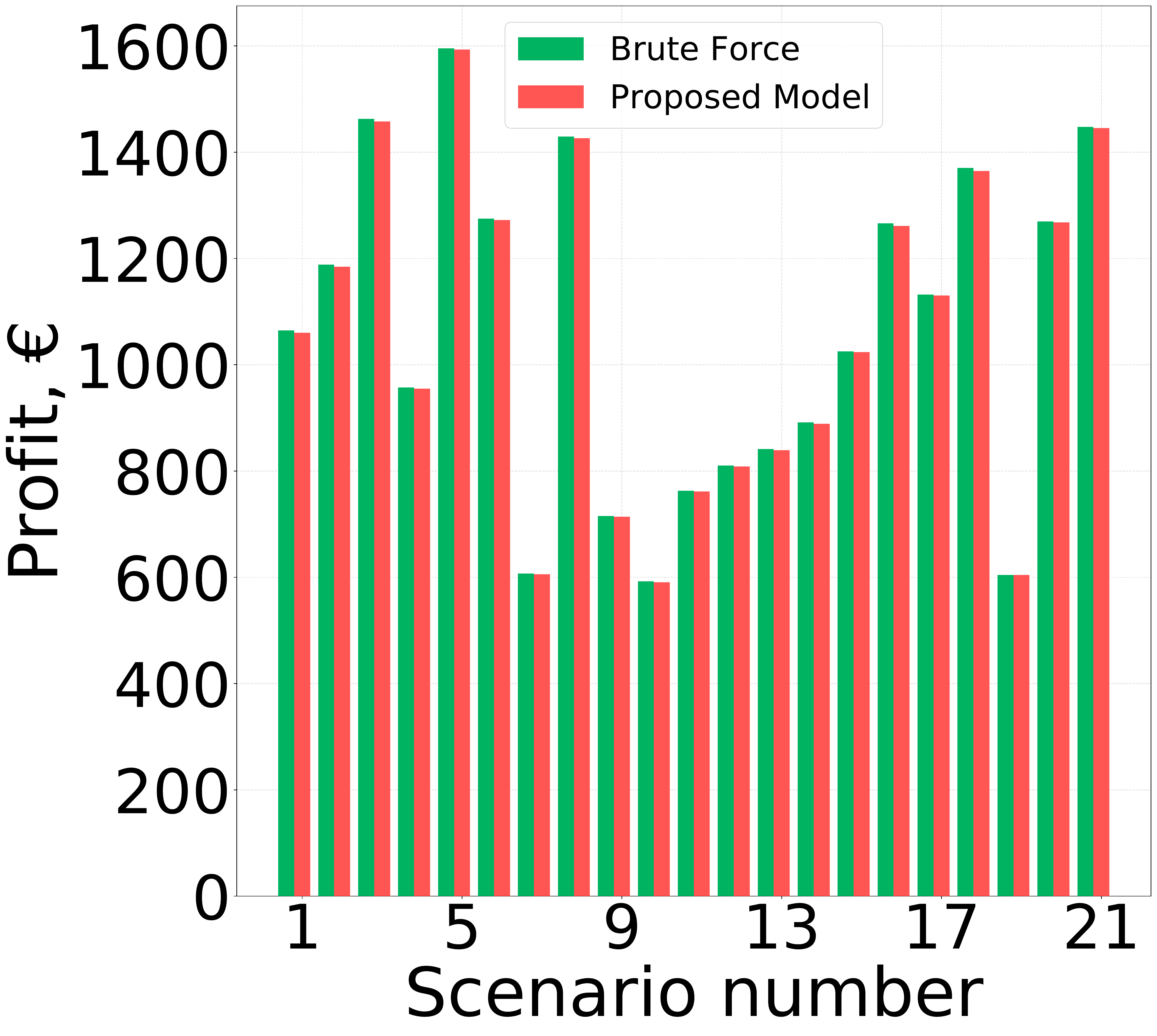} 
            } 
            \subfigure[Order placed for (\textit{Spring1, Realistic, Energy price 1}) scenario]
            {
                \label{subfig:Order validation}
                \includegraphics[width=0.23\textwidth]{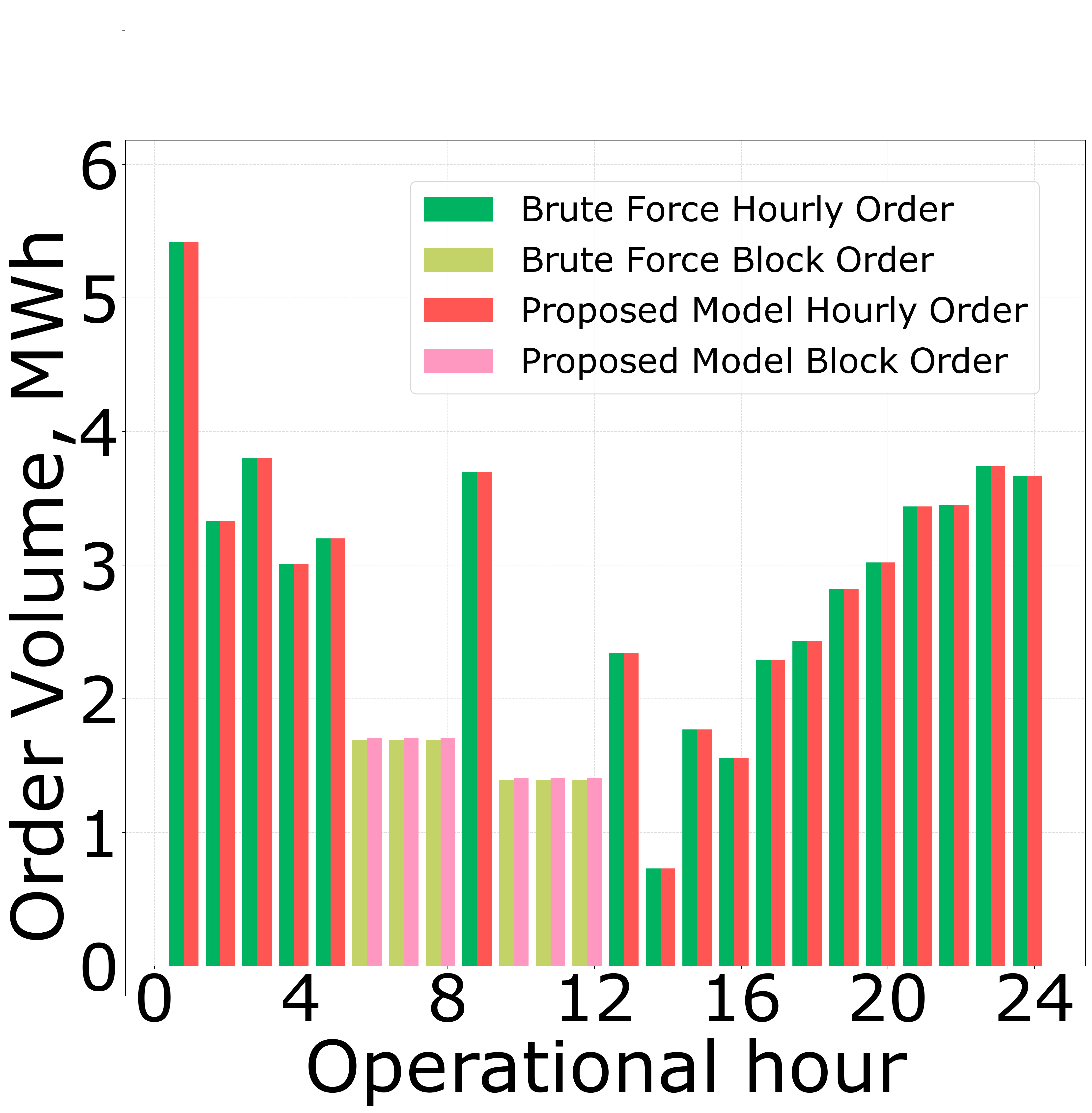} 
            } 
            \subfigure[Runtime comparison for acceptable scenarios]
            {
                \label{subfig:Runtime}
                \includegraphics[width=0.23\textwidth]{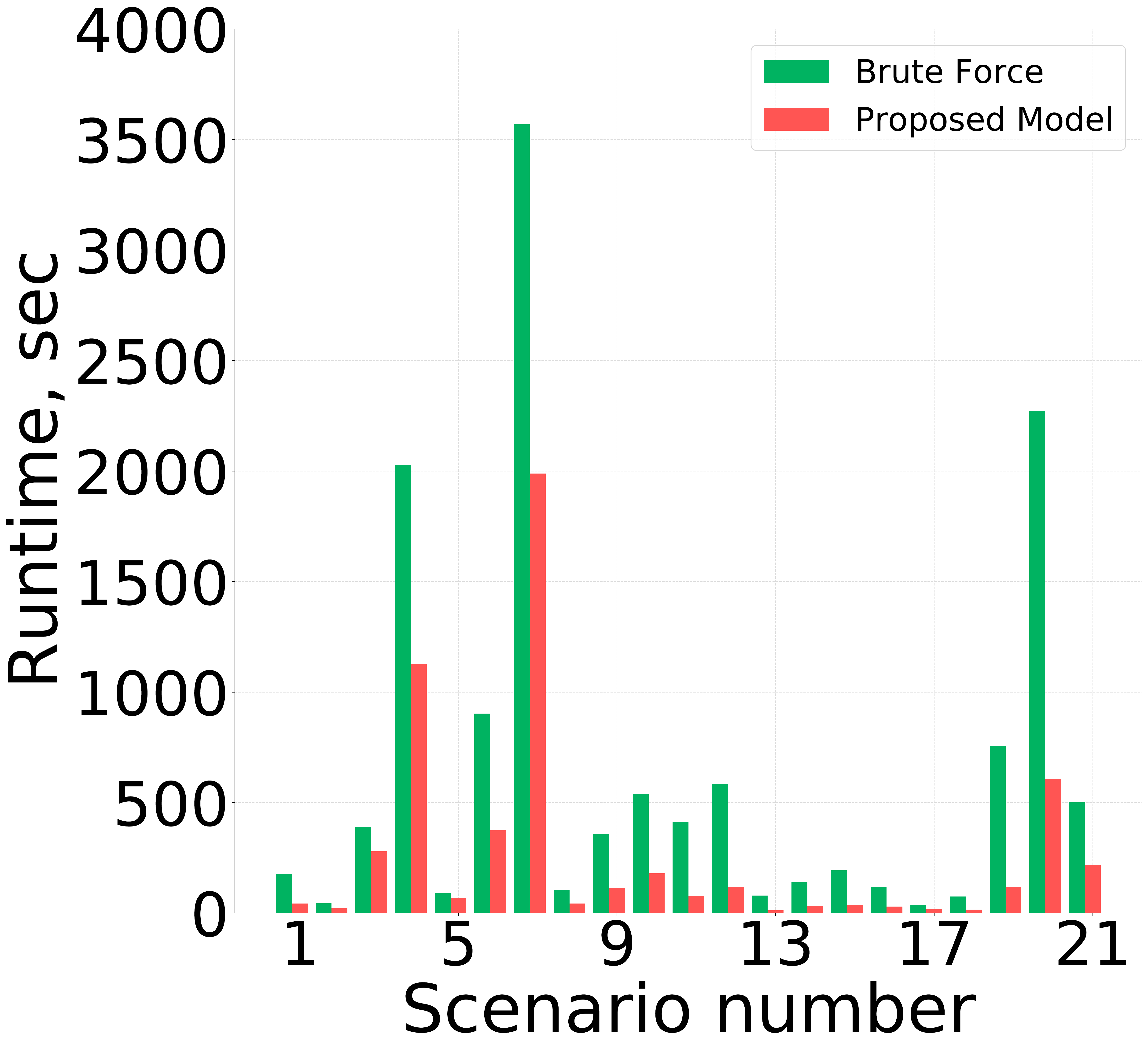}
            }
            \subfigure[\textit{MIP gap} vs. \textit{runtime} comparison for (\textit{Winter1, Realistic, Energy Price 1}) scenario] 
            {
                \label{subfig:MIP gap}
                \includegraphics[width=0.23\textwidth]{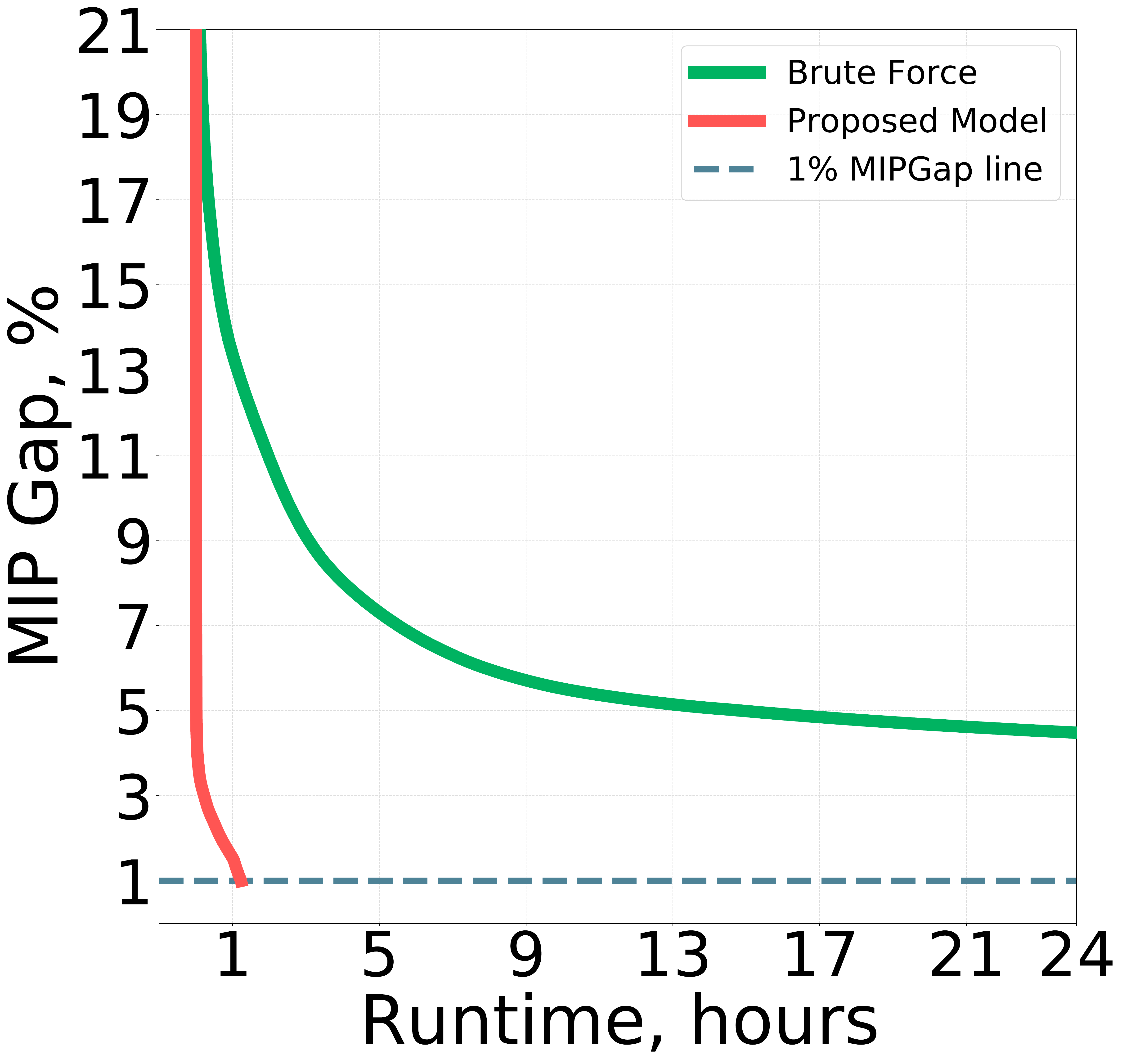}
            }
        \vspace{0.4em}
        \caption{Simulation results for different scenarios}
        \label{fig:Figure ref}
\end{figure*}

\section{Simulation Results}
\label{sec:Results}
The scenarios described in Section~\ref{sec:Simul} are simulated for the optimization formulations presented in Sections~\ref{subsec:Brute Force} and~\ref{subsec:AM2}. The simulations were run on a Windows machine with an Intel\textsuperscript{\textregistered} i7 8-core processor and 8 GB RAM using Python 3.8.8. Gurobi\textsuperscript{\textregistered} was selected as the solver for these simulations because it can solve to global optimum solutions for MIBLP and has a well-documented Python API~\cite{Gurobi21}. As mixed-integer programming (MIP) problems are NP-hard, two stopping criteria were used for the simulations conducted; a \textit{MIP gap} of \(1\%\) and a \textit{time limit} of \(3600\) sec, whichever occurs first. Thus, the solutions obtained for all scenarios are divided into two categories, namely acceptable (1\% \textit{MIP gap} within 3600 sec) or time cut-off (runtime exceeding 3600 sec) solutions. 
We observe that the proposed model could solve \textit{eight} more cases than brute force model within one hour. Also, all 21 common cases were acceptably (within 1\% \textit{MIP gap}) solved by brute force and proposed models. Thus, these 21 cases are used to make comparisons between the two models.

It is necessary to validate the performance of the models proposed in Section~\ref{subsec:AM2} using the brute force model from Section~\ref{subsec:Brute Force} to ensure the correctness of the proposed model. Figure~\ref{subfig:Obj comparison} shows the profit obtained (objective value) in the 21 common scenarios for each model. It can be seen that the objective values are about the same for all the cases. Additionally, Fig.~\ref{subfig:Order validation} shows the orders placed in the wholesale market by the two models for a particular scenario. It can be clearly seen that exactly the same orders are placed in the market by them, validating the proposed model. Additionally, for all the acceptable cases of each model, block orders were placed in the market. This shows that there exists opportunities for aggregators to benefit more by these types of orders, and modeling these is useful, especially in intervals, where thermal load flexibility is unavailable. 

In terms of computation speed, the proposed model outperforms the brute force formulation by 240\% with average runtimes of 263.30 and 637.08 sec, respectively. This was expected since MIP computational performance get exponentially worse with the increasing number of binary variables. Additionally, one of the \textit{time cut-off} cases was run for 24 hours for both methods and their \textit{MIP gaps} were compared in Fig.~\ref{subfig:MIP gap}. It can be clearly seen that our proposed approach reaches the 1\% \textit{MIP gap} in about 70 minutes whereas the brute force formulation is approaching 4.5\% \textit{MIP gap} after 24 hours, demonstrating the superiority of our proposed formulation. It also highlights that for certain scenarios, the brute force method may not be able to solve the problem within 24 hours (which is needed) for the DAM.  

For the following analyses, we use results from the proposed model as it solves more scenarios. The average profits obtained for the \textit{optimistic, realistic} and \textit{pessimistic} cluster are \texteuro 1166.81, \texteuro 1062.84, and \texteuro 922.60, respectively (per day). As expected, the aggregator makes more profit in the case of \textit{optimistic} cluster than the \textit{realistic} and \textit{pessimistic} clusters. The average profit considering all acceptable cases is \texteuro 1050.52. Comparing this method with sole single orders may be unfair since they would not correctly represent the rebound effect and temporal characteristics of the flexibility sources in question. Thus from these results, block orders show potential for accruing profit for aggregators.
\section{Conclusion and future work}
\label{sec:Conclusion}
In this paper, we presented a novel approach to optimally offer BTM flexibility in the Elspot market, considering hourly and regular block orders. We compared the proposed formulation against a pre-existing brute force model, while considering a generic flexibility model with inter-temporal dependencies. We generated 48 scenarios that consider seasonal effects and diurnal variations on flexibility, prosumer clusters with different price responsiveness and different spot prices. It was found that the proposed approach was correct and about 2.4 times faster than the brute force method. We also observed that this multi-product business model has potential to generate profit (\texteuro 1050.52 per day on average) for the aggregator. As part of the future work, we plan to incorporate other types of block orders and flexi-order in the proposed model and to analyze their potential benefit for the aggregator. To deal with variable flexibility and electricity prices, we plan to incorporate predictions and associated uncertainty of some factors influencing flexibility (e.g., weather and energy prices) in our formulations in the future. We also plan to making the formulation faster allowing us to incorporate more products/markets within it.

\newpage
\bibliographystyle{IEEEtran}
\bibliography{sample_base}
\end{document}